\documentclass{article}

\usepackage{arxiv}

\usepackage[utf8]{inputenc} 
\usepackage[T1]{fontenc}    
\usepackage{amsmath}
\usepackage{bm}
\usepackage{booktabs}       
\usepackage{amsfonts}       
\usepackage{amssymb}
\usepackage{nicefrac}       
\usepackage{microtype}      
\usepackage{soul}
\usepackage[table]{xcolor}
\definecolor{medium_blue}{rgb}{0, 0, 0.804}
\usepackage{hyperref}
\hypersetup{colorlinks=true,linkcolor=medium_blue, citecolor=medium_blue, urlcolor=medium_blue}
\usepackage{leftidx}
\usepackage{graphicx}
\usepackage{threeparttable}
\usepackage{multirow}
\usepackage{float}

\title{Cell-based Maximum Entropy Approximants for Three-Dimensional Domains: Application in Large Strain Elastodynamics using the Meshless Total Lagrangian Explicit Dynamics Method}

\author{
  Konstantinos A. ~Mountris\thanks{[mail] kmountris@unizar.es \quad [url] https://www.mountris.org} \\
  Arag\'on Institute of Engineering Research, IIS Arag\'on\\
  University of Zaragoza\\
  Spain, Zaragoza, ZGZ 50018 \\
  \texttt{kmountirs@unizar.es} \\
  \And
  George C. ~Bourantas \\
  Intelligent Systems for Medicine Laboratory\\
  The University of Western Australia\\
  Australia, Perth, WA 6009 \\
  \texttt{george.bourantas@uwa.edu.au} \\
  \AND
  Daniel ~Mill\'an \\
  CONICET \& Facultad de Ciencias Aplicadas a la Industria \\
  Universidad Nacional de Cuyo \\
  Argentina, San Rafael, M 5600\\
  \texttt{rdanielmillan@gmail.com} \\
  \And
  Grand R. ~Joldes \\
  Intelligent Systems for Medicine Laboratory\\
  The University of Western Australia\\
  Australia, Perth, WA 6009 \\
  \texttt{grand.joldes@uwa.edu.au} \\
  \AND
  Karol ~Miller \\
  Intelligent Systems for Medicine Laboratory\\
  The University of Western Australia\\ \vspace{1ex}
  Australia, Perth, WA 6009 \\
  Institute of Mechanics and Advanced Materials \\
  Cardiff School of Engineering \\
  Cardiff University \\
  UK, Cardiff, CF10 3AT \\
  \texttt{karol.miller@uwa.edu.au} \\
  \And
  Esther ~Pueyo \\
  Arag\'on Institute of Engineering Research, IIS Arag\'on, \\
  CIBER-BBN\\
  University of Zaragoza\\
  Spain, Zaragoza, ZGZ 50018 \\
  \texttt{epueyo@unizar.es} \\
  \AND
  Adam ~Wittek \\
  Intelligent Systems for Medicine Laboratory\\
  The University of Western Australia\\
  Australia, Perth, WA 6009 \\
  \texttt{adam.wittek@uwa.edu.au} \\
}

\begin{document}
\maketitle

\newpage

\begin{abstract}
We present the Cell-based Maximum Entropy (CME) approximants in $E^3$ space by constructing the smooth approximation distance function to polyhedral surfaces. CME is a meshfree approximation method combining the properties of the Maximum Entropy approximants and the compact support of element-based interpolants. The method is evaluated in problems of large strain elastodynamics for three-dimensional (3D) continua using the well-established Meshless Total Lagrangian Explicit Dynamics (MTLED) method. The accuracy and efficiency of the method is assessed in several numerical examples in terms of computational time, accuracy in boundary conditions imposition, and strain energy density error. Due to the smoothness of CME basis functions, the numerical stability in explicit time integration is preserved for large time step. The challenging task of essential boundary conditions imposition in non-interpolating meshless methods (e.g., Moving Least Squares) is eliminated in CME due to the weak Kronecker-delta property. The essential boundary conditions are imposed directly, similar to the Finite Element Method. CME is proven a valuable alternative to other meshless and element-based methods for large-scale elastodynamics in 3D. A naive implementation of the CME approximants in $E^3$ is available to download at \texttt{\url{https://www.mountris.org/software/mlab/cme}}.

\end{abstract}


\section{Introduction} \label{sec:introduction}
Meshfree Methods (MMs) \cite{Garg2018,Wittek2016} have been proposed as an alternative to the widely used Finite Element Method (FEM) for applications in solid mechanics, due to their attractive property to alleviate the mesh requirement for the approximation of an unknown field function. In mesh-based techniques, such as FEM, the approximation accuracy is deteriorated in problems that involve large deformations due to the mesh distortion. MMs are more suited for simulations involving large deformations. In MMs, the spatial domain is discretized by a set of nodes arbitrarily distributed without any interconnectivity. Since the introduction of Smoothed Particle Hydrodynamics (SPH) \cite{Lucy1977,Gingold1982}, MMs proliferated with developments such as the Meshless Collocation Methods (MCMs) \cite{Wen2007,Zhang2000,Bourantas2018}, the Meshless Local Petrov-Galerkin (MLPG) \cite{Atluri1998} and the Element-Free Galerkin (EFG) \cite{Belytschko1994} methods. A detailed overview of meshfree methods and their advantages can be found in \cite{Garg2018,Liu2005}.

In MCMs, the strong form solution of a PDE system is approximated using nodal collocation on the field nodes. Such methods demonstrate high efficiency, easy implementation, and direct Essential Boundary Conditions (EBC) imposition. However, MCMs often suffer from low accuracy when natural boundary conditions are involved due to the lower precision in the approximation of high-order derivatives \cite{Libre2008}. In weak form MMs, such as the MLPG and EFG methods, the order of derivation is reduced. Natural boundary conditions are satisfied in a straightforward manner. Weak form MMs have demonstrated high numerical stability and accuracy in large strain problems \cite{Hu2011,Cheng2016}.

The Meshless Total Lagrangian Explicit Dynamics (MTLED) \cite{Horton2010,Li2016,Joldes2019} is an efficient variant of the EFG method that was introduced to deal with large deformations. In MTLED, all calculations refer to the initial configuration of the continuum (total-Lagrangian formulation). MTLED implements a central difference scheme for explicit time integration, where basis functions and spatial derivatives are computed once, during the pre-processing stage. In addition, the explicit time integration demonstrates advantages over implicit integration schemes, such as (i) straightforward consideration of geometrical and material nonlinearity; (ii) easy implementation; (iii) suitability for massive parallelism \cite{Horton2010}. Originally, basis functions approximation in MTLED was performed using the Moving Least Squares (MLS) scheme  \cite{Horton2010}. An improved version of the MLS with advanced approximation capability named as the Modified MLS (MMLS) has been introduced in \cite{Joldes2015a}. The drawback of the approximation schemes of the MLS group is that the approximation of the field function does not possess the Kronecker-delta property and the imposition of Essential Boundary Conditions (EBC) requires special treatment. 

Methods that modify the weak form to impose EBC (e.g., Lagrange multipliers, penalty methods) are not suitable for explicit time integration. Joldes et al. \cite{Joldes2017} have proposed the EBC imposition in explicit meshless (EBCIEM) method for EBC imposition in MTLED. In EBCIEM, a displacement correction applies by distributing force on EBC nodes through Finite Element shape functions. EBCIEM applies exact imposition of EBC up to machine precision. However, the need for a Finite Element layer on the EBC boundary imposes restrictions regarding the applicability of the method to arbitrary geometries. A simplified version (SEBCIEM) has been also proposed \cite{Joldes2017}, where the distributed forces are lumped on the EBC nodes.

The EBC correction requirement can be alleviated if basis functions that possess the Kronecker-delta property are used. Approximation schemes of the Maximum Entropy (MaxEnt) group possess a weak Kronecker-delta property, which allows imposing EBC directly as in FEM. Different MaxEnt variations have been proposed up to date. Sukumar in \cite{Sukumar2004} used the maximization of information entropy to formulate meshfree interpolants on polygonal elements. Aroyo and Ortiz in \cite{Arroyo2006} used a Pareto compromise between locality of approximation and maximization of information entropy to create the Local Maximum Entropy (LME) approximants. LME approximants use generalized barycentric coordinates based on Jayne’s principle of maximum entropy \cite{Jaynes1957} and provide a seamless transition between non-local approximation and simplicial interpolation on a Delaunay triangulation (linear interpolants in the context of FEM).

LME approximants have already found a large number of applications \cite{Rosolen2013a,Rosolen2013b,Ullah2013,Gonzalez2010,Millan2013}. LME approximants are smooth and nonnegative approximants with local support that possess the weak Kronecker-delta property in the sense that basis functions on the boundary are not influenced by internal nodes. The LME support domain "locality" is controlled by a non-dimensional parameter $\gamma$. A variational approach to adapt $\gamma$ in the context of LME has been proposed in \cite{Rosolen2009}. Cell-based Maximum Entropy (CME) approximants have been proposed as an alternative to LME in the Galerkin framework \cite{Millan2015}. CME capitalizes on the background mesh connectivity only to define compact support for basis functions on $N$-rings ($N \ge 1$) of connected cells while the basis functions values are computed without using the connectivity information. Up to date, CME approximants have been implemented to approximate field functions in $E^2$.

The purpose of this work is to extend the CME approximants in $E^3$ and use them in the context of three-dimensional (3D) large-strain elastodynamic problems. Due to the previously described advantages of the explicit time integration over the implicit time integration for the solution of large-strain problems, we introduce the CME approximants in the MTLED method. While MTLED is an efficient method to address large deformation problems, the limitation of the special EBC treatment in MMLS can be restrictive in some scenarios. Our objective is to simplify the MTLED method and enable direct EBC imposition by replacing the MMLS with the CME. The rest of this paper is structured as follows. In Section \ref{sec:maxent}, we provide a compact and cogent formulation of the Cell-based Maximum Entropy approximation in $E^2$, based on the minimum relative entropy framework. In Section \ref{sec:cmeE3}, we extend the Cell-based Maximum Entropy approximation in $E^3$. In Section \ref{sec:cme_mtled}, we present benchmark examples, illustrating the accuracy and robustness of the proposed scheme in simple and complex geometries. Finally, in Section \ref{sec:conclusions}, we discuss our conclusions.

\section{Cell-based Maximum Entropy Approximants Formulation In \textit{E}\textsuperscript{2}} \label{sec:maxent}

\begin{figure}[htbp]
    \centering
    \includegraphics[width=0.65\textwidth]{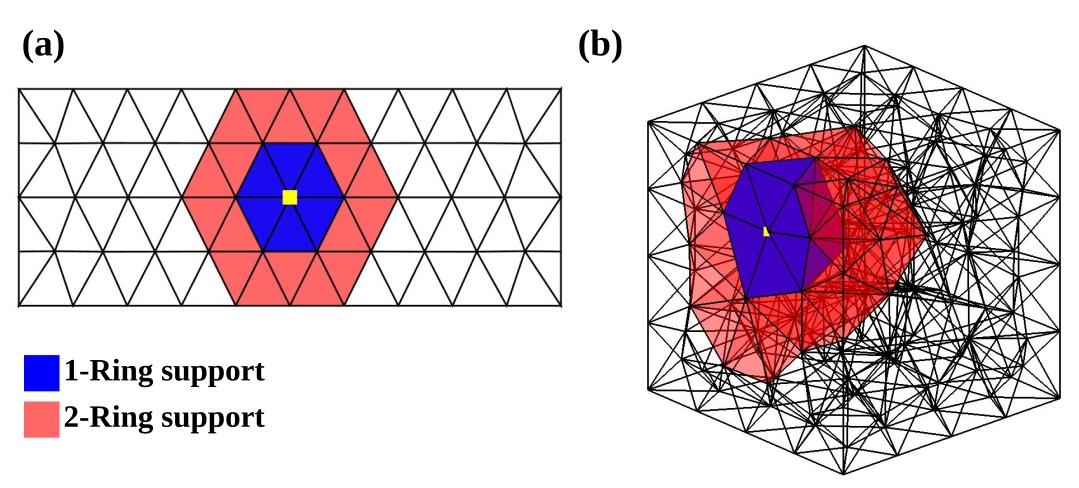}
    \caption{1\textsuperscript{st} and 2\textsuperscript{nd} ring support for a field node (yellow) in (a) two and (b) three dimensional domain.}
    \label{fig:cme_supports}
\end{figure}

CME approximants belong to the class of convex approximation schemes. They have compact support which leads to linear algebraic systems with small bandwidth and hence reduced memory requirements. CME are more time efficient compared to approximants with dense-nodal supports. Minimal supports are constructed for triangulated domains in $E^2$ by replacing Gaussian prior weights functions (commonly used in LME) with smooth distance approximation functions using the R-functions technology \cite{Shapiro1991}. The smooth distance approximation of a point $\bm{x} \in \Omega^2$ to the boundary $\partial \Omega^{N_R}$ of the $N_R$ support ring, $N_R = 1, 2,..., n \text{ and } \partial \Omega^{N_R} \setminus \partial \Omega$, is computed to derive the prior weight functions (Figure \ref{fig:cme_supports}).
Following, we give a brief overview of the minimum relative entropy framework (Subsection \ref{subsec:min_entropy}) along with the nodal prior weight construction using R-functions (Subsection \ref{subsec:cme_priors}). For a detailed description of the CME formulation we refer the reader to Mill\'an et al. \cite{Millan2015}.

\subsection{Minimum Relative Entropy Framework} \label{subsec:min_entropy}
For a scalar-valued function $u(\bm{x})$ and a set of unstructured points $\{\bm{x}_a\}_{a=1}^n \subset E^d$, the MaxEnt approximation $u^h(\bm{x})$ is given by

\begin{equation} \label{eq:maxent_approximation}
    u^h(\bm{x}) = \sum_{a=1}^n \phi_a (\bm{x}) u_a
\end{equation}

where $u_a$ are nodal coefficients, and $\{ \phi_a (\bm{x})\}_{a=1}^n$ are nonnegative basis functions that fulfill the zeroth and first order reproducing conditions:

\begin{equation} \label{eq:maxent_repconditions}
    \phi_a (\bm{x}) \geq 0, \quad \sum_{a=1}^n \phi_a (\bm{x}) = 1, \quad \sum_{a=1}^n \phi_a (\bm{x}) \bm{x}_a = \bm{x}
\end{equation}

The basis functions $\{ \phi_a (\bm{x})\}_{a=1}^n$ are defined from an optimization problem that is established at each evaluation point $\bm{x}$ for the linear constraints given by Equation (\ref{eq:maxent_repconditions}). Due to the nonnegative and partition of unity (zeroth order reproducing condition) properties, the basis functions can be interpreted as a discrete probability distribution. The informational entropy provides a canonical measure to the uncertainty associated with a discrete probability distribution. The least-biased approximation scheme that is consistent with the linear constraints is provided by the principle of maximum entropy. The formulation of MaxEnt approximants is derived by maximizing the Shannon-Jaynes entropy measure \cite{Jaynes1957} when nodal prior weights are used:

\begin{equation} \label{eq:shannon_entropy_measure}
    H ( \phi | w )= - \sum_{a=1}^n \phi_a (\bm{x}) \ln \bigg( \frac{\phi_a (\bm{x})}{w_a (\bm{x})} \bigg)
\end{equation}

where $D(\phi | w) = -H(\phi | w) \geq 0$ is the relative entropy measure and the corresponding variational principle is given by the principle of minimum relative (cross) entropy \cite{Shore1980}. The MaxEnt basis functions optimization problem can be stated in the variational formulation:

\begin{equation} \label{eq:maxent_dual}
\begin{aligned}
    (ME)_w \quad \text{For fixed } \bm{x} \text{, maximize } &H(\phi | w) \\
    \text{subject to } &\phi_a (\bm{x}) \geq 0  (a=1,\ldots,n), \\
                       &\sum_{a=1}^n \phi_a (\bm{x}) = 1, \\
                       &\sum_{a=1}^n \phi_a (\bm{x}) \bm{x}_a = \bm{x}
\end{aligned}
\end{equation}

\begin{figure}[htbp]
    \centering
    \includegraphics[width=\textwidth]{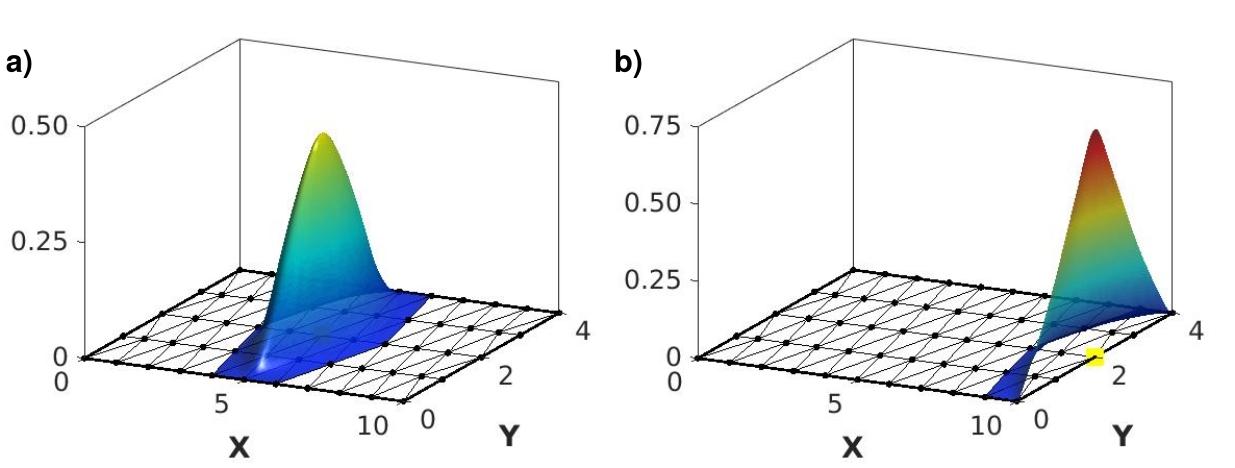}
    \caption{Cell-based Maximum Entropy (CME) basis function for a point located (a) at the center, with coordinates $\bm{x}=(5,2)$, and (b) at the right edge, with coordinates $\bm{x}=(10,2)$, of a two-dimensional rectangular domain.}
    \label{fig:cme_basis}
\end{figure}

Duality methods can be employed to solve the convex optimization problem in Equation (\ref{eq:maxent_dual}) efficiently and robustly \cite{Shore1980,Arroyo2007}. The basis functions are then:

\begin{equation} \label{eq:maxent_basis}
    \phi_a (\bm{x}) = \frac{w_a (\bm{x}) \text{ exp} [\bm{\lambda}(\bm{x}) \cdot (\bm{x} - \bm{x}_{a})]}{Z(\bm{x},\bm{\lambda}(\bm{x}))}
\end{equation}

where

\begin{equation} \label{eq:maxent_partition_func}
    Z(\bm{x},\bm{\lambda}(\bm{x})) = \sum_{b=1}^n w_b(\bm{x}) \text{ exp} [\bm{\lambda}(\bm{x}) \cdot (\bm{x} - \bm{x}_b)]
\end{equation}

is the partition function and $\bm{\lambda}$ is the Lagrange multiplier vector. The Lagrange multiplier vector is obtained solving the unconstrained convex optimization problem:

\begin{equation}
    \bm{\lambda}^* (x) = \text{arg min}_{\bm{\lambda} \in E^d} \ln Z(\bm{x},\bm{\lambda})  
\end{equation}

where $\bm{\lambda}^*$ is the converged solution of $\bm{\lambda}$ at $\bm{x}$. The resulting basis functions from the solution of $(ME)_w$ are noninterpolating, except on the boundary of the nodal set’s convex hull where the weak Kronecker-delta property holds (Figure \ref{fig:cme_basis}). The weak Kronecker-delta property provides a strong advantage over other approximant types, such as the MLS, where special treatments are required for the imposition of EBC \cite{Joldes2017}. Moreover, the basis functions inherit the nodal prior weight functions smoothness \cite{Sukumar2007,Arroyo2007}. Various nodal prior weight functions, such as Gaussian weight function \cite{Arroyo2006,Rosolen2013a,Li2010}, quartic polynomial weight function\cite{Yaw2009,Ortiz2010,Hale2012}, level set based nodal weight function \cite{Hormann2008,Sukumar2013}, exponential nodal weight function \cite{Nissen2012,Wu2014}, and approximate distance function to planar curves \cite{Millan2015} have been used to construct maximum entropy approximants with specific desired properties. The expression of the gradient for the MaxEnt basis functions for a node $a$ evaluated on a point $\bm{x} \in E^d$ is given by:

\begin{equation} \label{eq:maxent_grad}
    \nabla \phi_a^* = g_a^* \nabla w_a + \phi_a^* \left[ D \bm{\lambda}^* \cdot (\bm{x} - \bm{x}_a ) - \sum_{b=1}^n g_b^* \nabla w_b \right]
\end{equation}

where the superscript $*$ on any function indicates evaluation at $\bm{\lambda}^* (\bm{x})$, $w_a$ is the prior weight function for node $a$, $g_a^* := \phi_a^* / w_a$ , and $D \bm{\lambda}^*$ is given by:

\begin{equation}
    D \bm{\lambda}^* = - \bigg(\sum_{a=1}^n g_a^* \nabla w_a \otimes (\bm{x} - \bm{x}_a ) + \bm{I} \bigg) (\bm{J}^*)^{-1},
\end{equation}

where $\bm{I}$ is the identity matrix, $\bm{J}^* := \partial r / \partial \bm{\lambda}$, and $r := \sum_{a=1}^n \phi_a (\bm{x} - \bm{x}_a)$. A detailed description of the derivation of $\nabla \phi_a^*$ is given in Refs. \cite{Millan2015,Yaw2009}.

\subsection{Nodal Prior Weights Construction For CME In \textit{E}\textsuperscript{2}} \label{subsec:cme_priors}
To satisfy the required properties of CME (i.e., smoothness, compact support, unimodality), smooth approximations of the distance function to planar curves \cite{Biswas2004} are considered in the construction of prior weight functions. The theory of R-functions \cite{Shapiro1991} is used to approximate the distance function of each node in the nodal support of a basis function to the boundary polygon of the $N_R$-ring support domain (Figure \ref{fig:cme_supports}). Any real-valued function $F(\omega_1,\omega_2,\ldots,\omega_q)$, where $\omega_i (\bm{x}) : E^2 \rightarrow E$, is considered a R-function if its sign is determined only by the sign of its arguments $\omega_i$ and not their magnitude. Similarly to logical functions, R-functions can be written as a composition of elementary R-functions using operations such as negation, disjunction, conjunction, and equivalence \cite{Shapiro1991,Biswas2004}.

Using the composition property of R-functions, the approximate distance function of any given node $a$ from the boundary $\partial \Omega_a^{N_R}$ of its polygonal support domain $\Omega_a^{N_R}$ is expressed as the composition of elementary R-functions corresponding to the piecewise linear segments of the polygonal curve $\partial \Omega_a^{N_R}$. For each line segment belonging to $\partial \Omega_a^{N_R}$ with endpoints $\bm{x}_1 \equiv (x_1,y_1)$ and $\bm{x}_2 \equiv (x_2,y_2)$, the signed distance function from a point $\bm{x}$ to the line passing from $\bm{x}_1$ and $\bm{x}_2$ is defined by:

\begin{equation} \label{eq:rfun_sdf}
    f \equiv f(\bm{x}) := \frac{(x - x_1)(y_2 - y_1) - (y - y_1) (x_2 - x_1)}{L}
\end{equation}

where $L =  \parallel \bm{x}_2 - \bm{x}_1  \parallel $ is the length of the line segment. In addition, the line segment defines a disk of radius $L/2$. This disk can be expressed by the trimming function:

\begin{equation} \label{eq:rfun_trim}
    t \equiv t(\bm{x}) := \frac{1}{L} \left[\frac{L}{2}^2 - \parallel \bm{x} - \bm{x}_c \parallel \right]
\end{equation}

where $t$ is normalized to first order and $t \geq 0$ defines the disk with center $\bm{x}_c = (\bm{x}_1 + \bm{x}_2 )/2$. From Equations (\ref{eq:rfun_sdf}), (\ref{eq:rfun_trim}) the signed distance function from a point $\bm{x}$ to the line segment with endpoints $\bm{x}_1$ and $\bm{x}_2$ is given by the first order normalized function $\rho (\bm{x})$:

\begin{equation} \label{eq:rfun_normfun}
    \rho \equiv \rho(\bm{x}) :=  \sqrt{f^2 + \frac{1}{4} \bigg(  \sqrt{t^2 + f^4} - t \bigg)^2}
\end{equation}

where $\rho(\bm{x})$ is differentiable for any point that does not lay on the line segment. For $n$ line segments $l$ s.t. $\partial \Omega_a^{N_R} = l_1 \bigcup l_2 \bigcup \ldots \bigcup l_n$ the normalized approximation, up to order $m$, of the distance function $d_a$ for the node $a$ is given by the R-equivalence formula \cite{Biswas2004}:

\begin{equation} \label{eq:rfun_equivformula}
    d_a (l_1,l_2,\ldots,l_n ) := \rho_1 \sim \rho_2 \sim \ldots \sim \rho_n = \frac{1}{\sqrt[m]{\frac{1}{\rho_1^m} + \frac{1}{\rho_2^m} + \ldots + \frac{1}{\rho_n^m}}}
\end{equation}

The R-equivalence formula has the desirable properties of preserving the normalization up to the $m^{th}$  order and being associative, while the R-conjunction is normalized up to the $(m-1)^{th}$ order and it is not associative. Finally, the prior weight function $w_a(\bm{x})$ for node $a$ is obtained by:

\begin{equation} \label{eq:cme_prior}
    w_a(\bm{x}) = \frac{d_a^s (\bm{x})}{\sum_{b \in N_x} d_b^s(\bm{x})}
\end{equation}

where $N_x$ are the indices of the $N_R$-ring nodal neighbors of the point $\bm{x}$, $s \geq 2$ is a smoothness modulation factor and $w_a(\bm{x})$ is normalized s.t. $0 \leq w_a(\bm{x}) \leq 1$. It should be noted that the R-equivalence formula in Equation (\ref{eq:rfun_equivformula}) does not preserve the normalization of the distance function on the joining vertices of the piecewise linear segments and hence introduces undesired bulging effect at these points. Increasing the smoothness modulation factor $s$, the bulging effect is reduced.

\section{Extension Of Cell-based Maximum Entropy Approximants In \textit{E}\textsuperscript{3}} \label{sec:cmeE3}

MaxEnt approximants (Equation (\ref{eq:maxent_basis})) are naturally generalized to $E^d$, with $d \geq 3$. However, CME approximants have been limited to $E^2$ due to the construction of nodal prior weight functions as smooth approximations to the distance function of 2D polygonal curves. In this section we extend the CME approximants to $E^3$ by constructing nodal prior weight functions as smooth approximations to the distance functions of polyhedral surfaces following the development of Biswas and Shapiro \cite{Biswas2004}. 

\subsection{Nodal Prior Weights Construction For CME In \textit{E}\textsuperscript{3}} \label{subsec:cme_priors_E3}

We consider the discretization $\Omega^h$ of a finite domain $\Omega \in E^3$ in tetrahedra $\tau$. The $N_R$  ring support domain for any vertex $a$ of $\tau \in \Omega^h$ is defined by the union of the piecewise triangular faces belonging to the $N_R$ ring of attached tetrahedra to the vertex $a$ (Figure \ref{fig:cme_supports}). R-functions are employed to construct the approximation distance field to the boundary of the $N_R$ ring support domain by joining the fields of the individual boundary triangular faces.

\begin{figure}[htbp]
    \centering
    \includegraphics[width=0.65\textwidth]{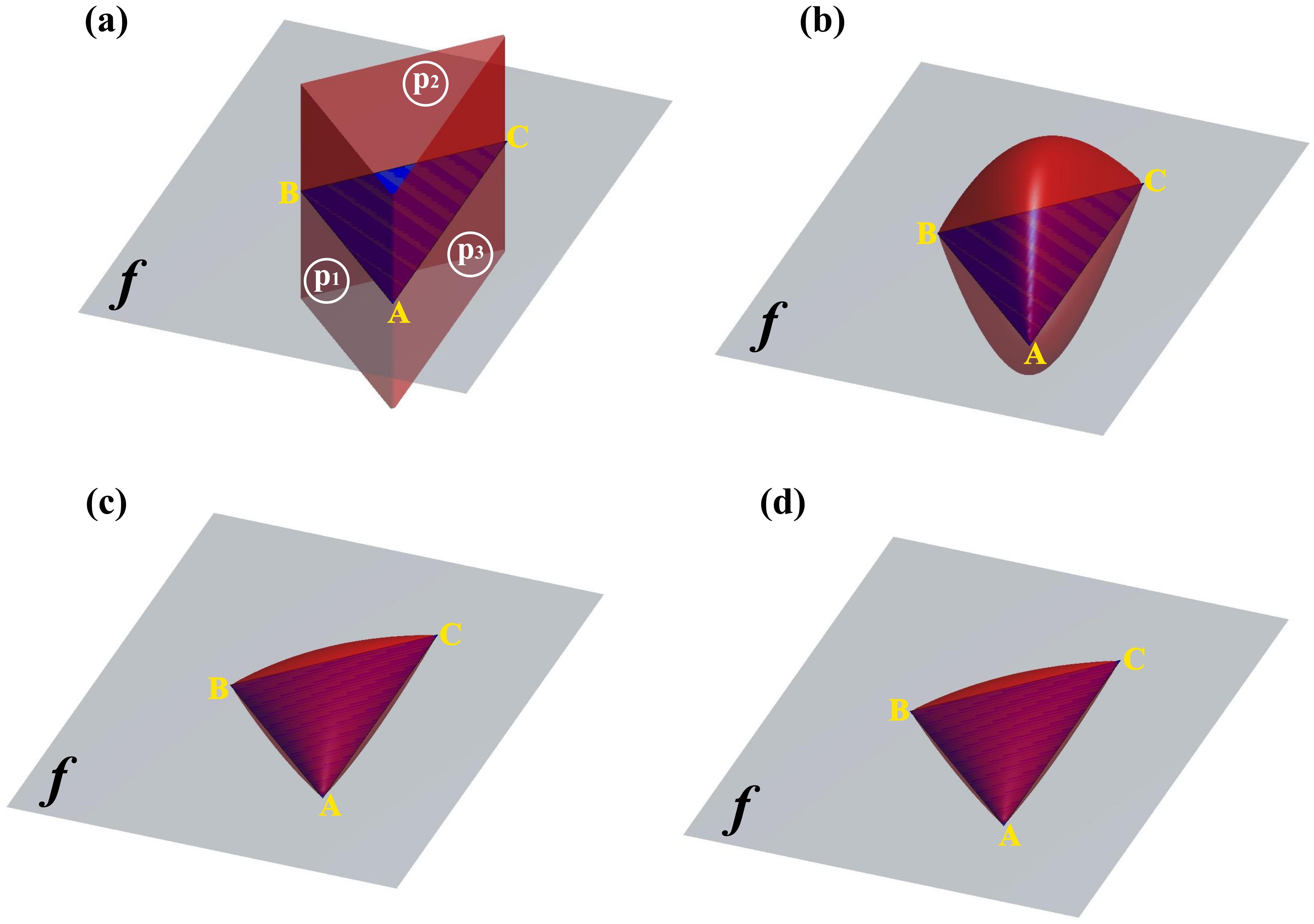}
    \caption{(a) Prism trim volume for a triangle constructed by planes $p_1$, $p_2$, $p_3$ ($a = 0$). (b - d) Smooth trimming volume tangent to the edges of the triangle for $a = \{0.2, 1.6, 2.0\}$. Parameter $k=2$.}
    \label{fig:cme_E3_trim}
\end{figure}

The distance field approximation to any triangle is constructed considering the intersection of the triangle’s carrier plane $f(\bm{x}) = 0$ and a trim volume $t$ defined by the planes $p(\bm{x})_i$, $i=1,2,3$ that pass through the edges of the triangle and are orthogonal to the carrier plane $f(\bm{x})$ (Figure \ref{fig:cme_E3_trim}). The trim volume $t$ can be constructed joining the planes $p_i$ according to the $R_f$-conjunction formula:

\begin{equation} \label{eq:rfun_conjunc}
    p_1 \wedge_f p_2 \equiv p_1 + p_2 - \sqrt{p_1^2 + p_2^2 + \alpha f^k}
\end{equation}

where $\alpha$, $k$ are constants controlling the shape of the R-function’s zero set. For points in $E^3$ that do not lay on the carrier plane the Equation (\ref{eq:rfun_conjunc}) is $C^\infty$, while for points on the carrier plane ($f=0$) it reduces to the standard $R$-conjunction formula. Therefore, $t$ is given by:

\begin{equation} \label{eq:rfun_trimE3}
    t \equiv (p_1 \wedge_f p_2) \wedge_f p_3 \equiv p_1 + p_2 + p_3 - G - \sqrt{(p_1 + p_2 - G)^2 + p_3^2 + \alpha f^k}
\end{equation}

where $G := \sqrt{p_1^2 + p_2^2 + \alpha f^k}$. The function $t$ defines a smooth trimming volume tangent to the three edges of the triangle. The coefficient $\alpha$ controls the smoothness of the trimming volume, where for large $\alpha$ the trimming volume flattens out and for $\alpha = 0$ it becomes identical to the unbounded prism constructed by the combination of the planes $p_i$ (Figure \ref{fig:cme_E3_trim}). While the unbounded prism trimming volume ($\alpha = 0$) does not maintain differentiability on the triangle edges, the smooth volume ($\alpha > 0$) is differentiable at any point except the triangle vertices where the bulging effect arises. Substituting $f$ and $t$ in Equation (\ref{eq:rfun_normfun}), the normalized distance function from any point $\bm{x}$ to the triangle is approximated and the nodal prior weights $w_a (\bm{x})$ are obtained by Equations (\ref{eq:rfun_equivformula}), (\ref{eq:cme_prior}), similarly to the $E^2$ case.

\subsection{Derivation Of The CME Approximants Gradient In \textit{E}\textsuperscript{3}} \label{subsec:cme_grad_E3}

To compute the gradient of the CME approximants $\nabla \phi_a (\bm{x})$ (Equation (\ref{eq:maxent_grad})), we should first compute the gradient of the prior function $\nabla w_a (\bm{x})$. Following the abbreviation in \cite{Millan2015}, to derive $\nabla w_a (\bm{x})$, for a node $a$  at a point $\bm{x}$, we first derive the gradient of the normalized approximation to the distance function:

\begin{equation} \label{eq:ndf_grad}
    \nabla d_a(\bm{x}) \equiv \nabla d_a = \frac{\sum_{i=1}^{n_a} \rho_i^{-(m+1)} \nabla \rho_i}{(\sum_{i=1}^{n_a} \rho_i^{-m})^{\frac{m+1}{m}}} 
\end{equation}

where $\rho_i$ is the distance field function of the $i^{th}$ triangular patch of the support domain’s boundary for node $a$. The gradient of the distance field function $\rho$ is given by:

\begin{equation} \label{eq:rho_grad}
    \nabla \rho = \frac{1}{\rho} \Bigg(f \nabla f + \left[ \frac{\sqrt{t^2+f^4}-t}{4} \right] \left[\frac{t \nabla t + 2f^3  \nabla f}{\sqrt{t^2 + f^4}} \right] \Bigg)
\end{equation}

where $i$ is omitted for notation simplicity, $\nabla f = \widehat{n_f} = A_f \widehat{i} + B_f \widehat{j} + C_f \widehat{k}$ is the gradient of the triangle’s carrier plane $f$, which is identical to the plane’s normal vector, and $\nabla t$ is the gradient of the trim volume function:

\begin{equation} \label{eq:trimfun_grad}
\begin{aligned}
    \nabla t &= \nabla \bigg(p_1 + p_2 + p_3 - G - \sqrt{(p_1 + p_2 - G)^2 + p_3^2 + \alpha f^k} \bigg) \\
    \nabla t &= \nabla(p_1 + p_2 + p_3 - G) - \nabla \bigg(\sqrt{(p_1 + p_2 - G)^2 + p_3^2 + \alpha f^k} \bigg) \\
    \nabla t &= \nabla (p_1 + p_2 + p_3 - G) - \frac{\nabla (p_1 + p_2 - G)^2 + \nabla p_3^2 + \alpha \nabla f^k}{2\sqrt{(p_1 + p_2 - G)^2 + p_3^2 + \alpha f^k}} \\
    \nabla t &= \nabla p_1 + \nabla p_2 + \nabla p_3 - \nabla G - \frac{2(p_1 + p_2 - G)(\nabla p_1 + \nabla p_2 - \nabla G) + 2p_3 \nabla p_3 + \alpha k f^{k-1} \nabla f}{2\sqrt{(p_1 + p_2 - G)^2 + p_3^2 + \alpha f^k}}
\end{aligned}
\end{equation}

where $\nabla p_j$, $j = 1,2,3$ is the gradient of the perpendicular plane to the $j^{th}$ edge of the triangle given by its normal vector $\widehat{n_{p_j}}$ and $\nabla G$ is given by:

\begin{equation} \label{eq:G_grad}
\begin{aligned}
    \nabla G &= \nabla \bigg( \sqrt{p_1^2 + p_2^2 + \alpha f^k} \bigg) \\
    \nabla G &= \frac{2p_1 \nabla p_1 + 2p_2 \nabla p_2 + \alpha k f^{k-1} \nabla f}{2\sqrt{p_1^2 + p_2^2 + \alpha f^k}}
\end{aligned}
\end{equation}

\section{Evaluation Of The CME Approximants In The MTLED Method} \label{sec:cme_mtled}

We construct the matrix of the basis function derivatives using the CME approximation, instead of the MMLS, to derive the CME-MTLED method. We perform a series of numerical examples for three-dimensional large strain hyperelasticity at steady state to assess the efficiency of the CME approximants in the MTLED method. In all the examples, we use the hyper-elastic neo-Hookean material which is described by the hyper-elastic strain energy density function \cite{Bonet2008}:

\begin{equation} \label{eq:neohookean_law}
    W = \frac{\mu}{2}(I_1-3) - \mu \ln J + \frac{\lambda}{2} (\ln J)^2
\end{equation}

where $\lambda$ and $\mu$ are the Lam\'e parameters, $I_1$ is the first strain invariant and $J$ is the determinant of the deformation gradient. The evaluation of spatial integrals is applied on quadrature points, generated on a background unstructured tetrahedral mesh using a 4-point quadrature rule \cite{Cools1993}. We use dynamic relaxation, described in detail in \cite{Joldes2009,Joldes2011}, to ensure fast convergence for the explicit time integration scheme to the steady state solution. All the simulations are performed using a multi-threaded implementation of the MTLED algorithm on an Intel Core i7-8700K CPU using 12 threads and 16GB DDR4 RAM.

\subsection{Cube Under Unconstrained Compression} \label{subsec:unconstrained_compression}

\begin{figure}[H]
    \centering
    \includegraphics[width=0.5\textwidth]{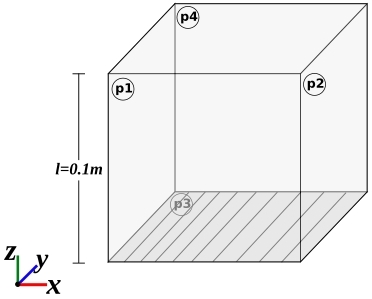}
    \caption{Boundary conditions for unconstrained cube under compression. Face \textbf{p1}: $u_y = 0$, $u_x$ and $u_z$ unconstrained, face \textbf{p2}: $u_x = 0$, $u_y$ and $u_z$ unconstrained, face \textbf{p3}: $u_z = 0$, $u_x$ and $u_y$ unconstrained, and face \textbf{p4}: $u_z = -0.02 m$, $u_x$ and $u_y$ unconstrained.}
    \label{fig:cube_unconstrained}
\end{figure}

We model the compression of an unconstrained cube with soft-tissue-like constitutive properties (hyper-elastic neo-Hookean material). The cube has side length $l=0.1 m$, Young modulus $YM =3000 Pa$, Poisson ratio $v=0.49$, and density $\rho = 1000 kg/m^3$. It is subjected to unconstrained compression (Figure \ref{fig:cube_unconstrained}) by displacing the top surface by $0.02 m$ ($20\%$ of the initial height). In this simple case, the vertical component of displacement ($u_z$) is given by the analytical solution $u_z=-0.2z$, ($z$ refers to the reference configuration).

Simulations are performed for two successively fined nodal distributions, consisting of $152$ and $4594$ nodes, respectively. Computations are conducted using the CME-MTLED method. We evaluate the accuracy of the proposed scheme using the Normalized Root Mean Square Error ($NRMSE$):

\begin{equation} \label{eq:nrmse}
    NRMSE = \frac{\frac{1}{N}\sum_{i=1}^N (u_{z_i}^h - u_{z_i}^{an})^2}  {u_{z_{max}}^{an} - u_{z_{min}}^{an}}
\end{equation}

where $N$ is the number of nodes, $u_{z_i}^h$ is the $z$-axis displacement component of the numerical solution, and $u_{z_i}^{an}$ the $z$-axis displacement component of the analytical solution. Table \ref{tab:cube_unconst_errors} shows the computed $NRMSE$ for the two considered nodal distributions. Good agreement with the analytical solution is achieved even for the coarse grid ($152$ nodes, $1896$ quadrature points). Additionally, due to the weak Kronecker-delta property of the CME approximants, the essential boundary conditions defined on faces $p1$, $p2$, $p3$ and $p4$ are imposed exactly (absolute error is zero up to machine precision).

\begin{table}[H]
    \caption{Normalized Root Mean Square Error ($NRMSE$) for a cube under unconstrained compression.}
    \centering
    \begin{threeparttable}
        \begin{tabular}{c c c}
            \toprule
            Nodes  & Quadrature points & $NRMSE$ \\
            \midrule
            152    & 1896              & 1.75E-3 \\
            4594   & 92876             & 3.85E-4 \\
            \bottomrule
        \end{tabular}
    \end{threeparttable}
    \label{tab:cube_unconst_errors}
\end{table}

\subsection{Cube Under Constrained Deformation} \label{subsec:uniaxial_deformation}

\begin{figure}[htbp]
    \centering
    \includegraphics[width=\textwidth]{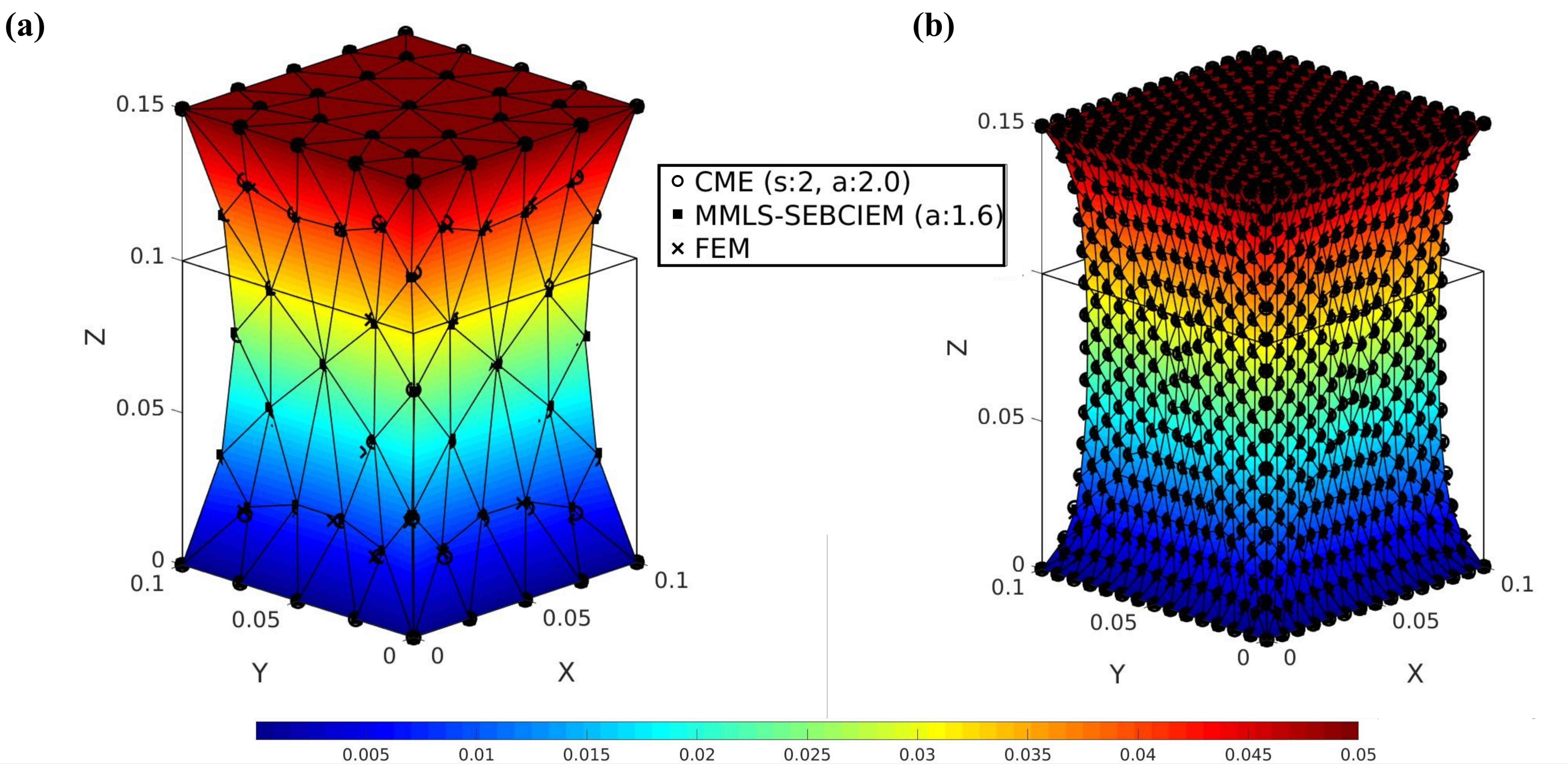}
    \caption{Steady state deformation under $50\%$ extension for (a) $lvl_1$, and (b) $lvl_4$ discretizations. Nodes in the bottom surface are fixed and nodes at the top surface undergo a displacement in the $z$-axis, $u_z = -0.05m$ $\big(\bigcirc$ CME, $\blacksquare$ MMLS-SEBCIEM, $\bm{\times}$ FEM$\big)$.}
    \label{fig:cube_extension}
\end{figure}

We further demonstrate the accuracy of the proposed scheme by considering the extension and compression of a cube with side length $l=0.1 m$. One face ($z=0 m$) of the cube is rigidly constrained, while the opposite face ($z=0.1 m$) is displaced (Figure \ref{fig:cube_extension}). A maximum displacement loading of $0.05 m$ ($50\%$ of the initial height) is smoothly applied. A hyper-elastic neo-Hookean material model with $YM=3000 Pa$, $v=0.49$ and $\rho=1000 kg/m^3$ is chosen to capture the material response of soft tissue.

The cube is discretized using six successively denser nodal distributions ($lvl_1$: $152$, $lvl_6$: $50521$). The CME approximants are used to approximate the unknown field functions. Simulations are performed to test the convergence of the proposed scheme. The simulation characteristics such as the critical explicit integration time step ($\Delta t_{crit}$), the number of execution steps ($N_{exe}$), and the total execution time ($t_{exe}$) are evaluated for CME with $N_R = 2$ and parameters $s=\{2,4\}$; $m=3$; $a=\{0.2,1.6,2.0\}$; $k=2$. The simulations are also performed using the MMLS approximants with the SEBCIEM correction \cite{Joldes2017}. The dilatation coefficient $a=1.6$ (parameter controlling support size) is used in MMLS. The simulation characteristics are given in Table \ref{tab:cube_uniaxial_results}.

The $\Delta t_{crit}$ for CME is found higher for all the parameter sets compared to MMLS-SEBCIEM. The $\Delta t_{crit}$ gain is reduced for increasing $s$ and $a$; with increasing $s$ the gradient becomes steeper, leading to larger eigenvalues in the spatial derivatives matrix ($\leftidx{^t_0}{\bm{B}}_{L0}^t$). The CME performs better, in terms of computational efficiency, since it can reach steady state with less time steps than the MMLS-SEBCIEM.

\begin{table}[htbp]
    \caption{Simulation characteristics (critical explicit integration time step - $\Delta t_{crit}$, number of execution steps - $N_{exe}$, execution time - $t_{exe}$) for cube under $50\%$ uniaxial extension $\&$ compression for $lvl_1$-$lvl_6$ discretization.}
    \centering
    \begin{threeparttable}
        \begin{tabular}{c c c c c c}
            \toprule
            Approx.  & $s$ \tnote{$\dagger$} & $a$ \tnote{$\ddagger$}  & $\Delta t_{crit}$ (ms) & $N_{exe}$ & $t_{exe}$ (s) \\
            \midrule
            \multicolumn{6}{c}{Extension} \\
            \midrule
            \multirow{6}{*}{CME}    &  \multirow{3}{*}{2}  & 0.2  & 2.173 - 0.260    & 2963 - 10455   & 1.4 - 3173.5  \\
            &  & 1.6  &  1.972 - 0.243  & 2814 - 12707  & 1.6 – 3906.9  \\
            &  & 2.0  &  1.918 - 0.238  & 3341- 13122  & 1.4 - 3911.0  \\
            &  \multirow{3}{*}{4}  & 0.2  & 1.527 - 0.205 & 3396 - 16582   & 1.9 - 4652.4  \\
            &  & 1.6  & 1.230 - 0.176  & 3529 - 20900  & 1.8 - 5711.7  \\
            &  & 2.0  & 1.216 - 0.174  & 3509 - 21167  & 1.9 - 6070.6  \\
            MMLS-SEBCIEM  &  1.6  & 0.2  & 0.902 - 0.192   & 3208 - 16827  & 1.4 - 4603.5  \\
            \midrule
            \multicolumn{6}{c}{Compression} \\
            \midrule
            \multirow{6}{*}{CME}    &  \multirow{3}{*}{2}  & 0.2 & 2.173 - 0.260  & 3172 - 10442   & 1.7 - 3184.9  \\
            &  & 1.6  & 1.972 - 0.243  & 3093 - 17288  & 1.5 - 4630.3  \\
            &  & 2.0  & 1.918 - 0.238  & 3117 - 16501 & 1.5 - 4414.7  \\
            &  \multirow{3}{*}{4}  & 0.2  & 1.527 - 0.205 & 3279 - 31402   & 1.7 - 7231.9  \\
            &  & 1.6  & 1.230 - 0.176  & 3664 - 24885 & 2.0 - 6297.4 \\
            &  & 2.0  & 1.216 - 0.174  & 3735 - 18563  & 1.7 - 5557.1 \\
            MMLS-SEBCIEM  &  1.6  & 0.2  & 0.902 - 0.192  & 3975 - 18315  & 1.8 - 5710.6  \\
            \bottomrule
        \end{tabular}
        \begin{tablenotes}
            \item [$\dagger$] CME smoothness modulation factor
            \item [$\ddagger$] CME R-function zero set shape factor or MMLS dilatation coefficient
        \end{tablenotes}
    \end{threeparttable}
    \label{tab:cube_uniaxial_results}
\end{table}

\begin{figure}[htbp]
    \centering
    \includegraphics[width=\textwidth]{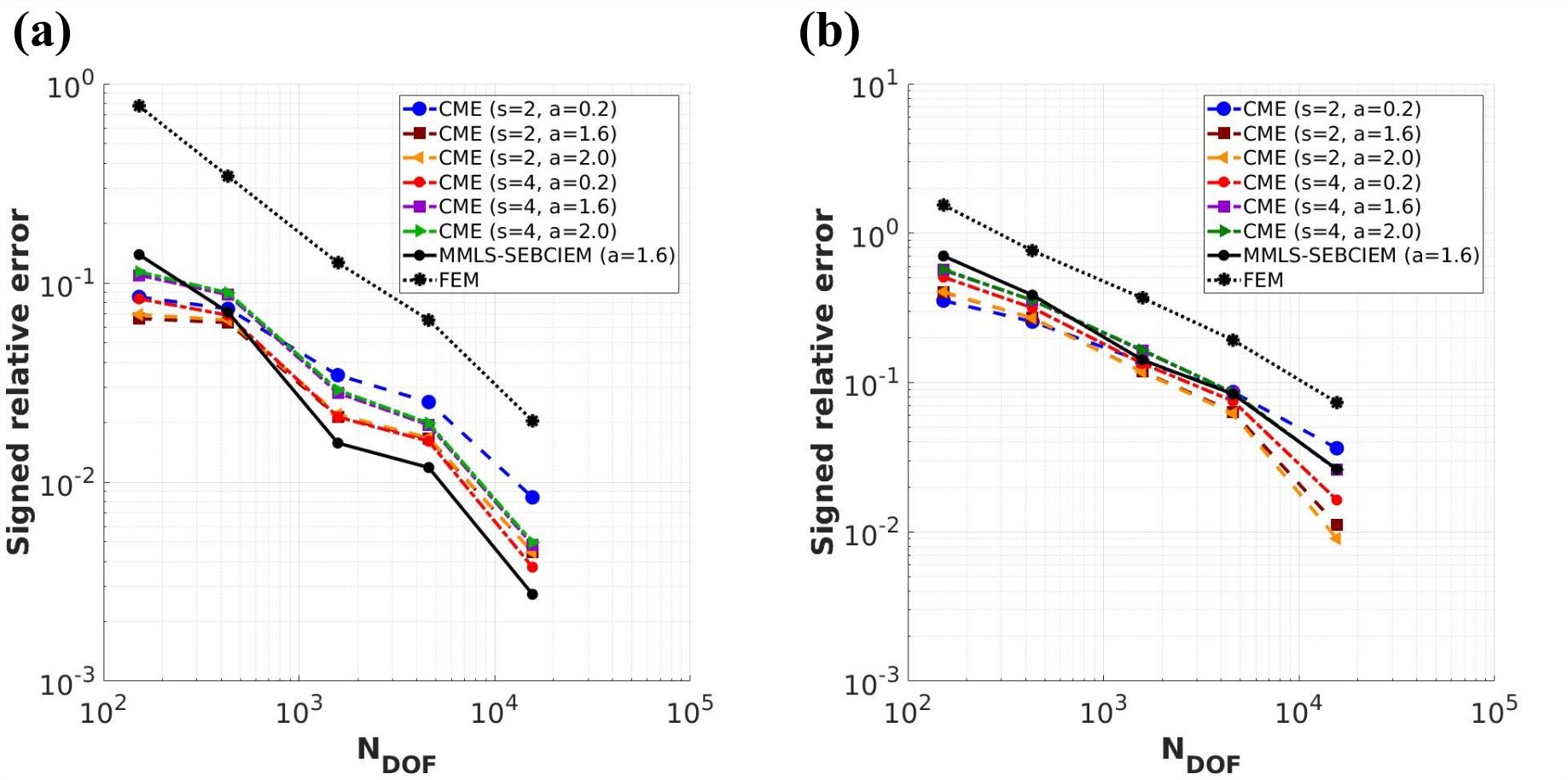}
    \caption{Signed relative error in strain energy density $SRE_W$ for a cube at: (a) $50\%$ uniaxial extension, (b) $50\%$ uniaxial compression.}
    \label{fig:cube_uniaxial_convergence}
\end{figure}

An analytical solution is not available for the given problem. Therefore, to evaluate the convergence and accuracy of the method we use the Signed Relative Error in strain energy density ($SRE_W$) as a convergence measure similar to \cite{Arroyo2006}. We compute the $SRE_W$ for $lvl_1$ to $lvl_5$ discretization with respect to the $lvl_6$ (finest) discretization (Figure \ref{fig:cube_uniaxial_convergence}). The $SRE_W$ is obtained by:

\begin{equation}
    SRE_W = \frac{\Bar{W}_{lvl_i} - \Bar{W}_{lvl_6}}{\Bar{W}_{lvl_6}}
\end{equation}

where $\Bar{W}_{lvl_i}$ is the mean value of the strain energy density for the $lvl_i$ discretization, $i=1,\ldots,5$. The $SRE_W$ for CME is comparable with the $SRE_W$ for MMLS-SEBCIEM. While the MMLS-SEBCIEM lead to lower $SRE_W$ for higher nodal density ($lvl_3$-$lvl_5$), the CME results in lower $SRE_W$ for lower nodal density ($lvl_1$ and $lvl_2$).

For additional comparison, the $SRE_W$ for $lvl_1$ to $lvl_5$ discretization for a Finite Element Method (FEM) simulation is given. FEM simulations are performed with the FEBio v$2.5.2$ software \cite{Maas2012} using isoparametric tetrahedral elements and the Newton-Raphson implicit time integration method. Linear tetrahedral elements are known for being prone to volumetric “locking”. For this reason, FEBio implements a nodally integrated tetrahedron with enhanced performance for finite deformation and near-incompressibility compared to the standard constant strain tetrahedron \cite{Puso2006}. The $SRE_W$ for FEM simulations is found an order of magnitude higher from CME in $50\%$ extension.

\begin{figure}[htbp]
    \centering
    \includegraphics[width=\textwidth]{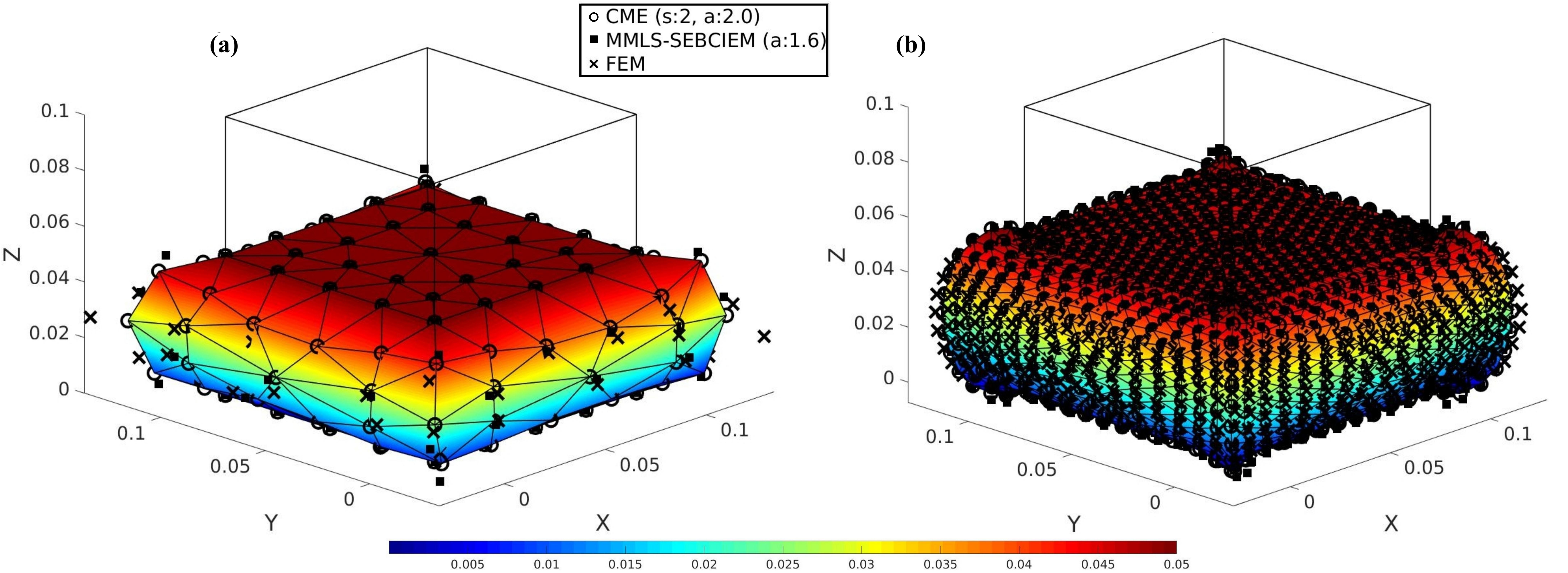}
    \caption{Steady state deformation under $50\%$ compression for (a) $lvl_1$, and (b) $lvl_4$ discretizations. Nodes in the bottom surface are fixed and nodes at the top surface undergo a displacement in the $z$-axis, $u_z = 0.05m$ $\big(\bigcirc$ CME, $\blacksquare$ MMLS-SEBCIEM, $\bm{\times}$ FEM$\big)$.}
    \label{fig:cube_compression}
\end{figure}

Similar results are acquired for $50\%$ compression. The use of CME in MTLED leads to higher $\Delta t_{crit}$ and lower $t_{exe}$ compared to the MMLS-SEBCIEM (Table \ref{tab:cube_uniaxial_results}). The $SRE_W$ in $50\%$ compression is found lower for CME compared to MMLS-SEBCIEM and FEM for all the discretization levels (Figure \ref{fig:cube_uniaxial_convergence}). The set-up of $N_R = 2$ and $s=2$; $m=3$; $a=2.0$; $k=2$ for CME is found to be a good trade-off between accuracy and efficiency in MTLED for both cases of extension and compression (Table \ref{tab:cube_uniaxial_results}). The use of CME instead of the MMLS-SEBCIEM in MTLED eliminates the necessity for EBC correction and the execution time is reduced. The execution time reduction is more evident for dense nodal discretization (up to $1295.9 s$ - $22.7\%$). Moreover, the evaluation of the $SRE_W$ demonstrates the improved accuracy of the MTLED method over the FEM for large strain problems for both $50\%$ extension and compression conditions. Especially for severe distortion, such as in the $50\%$ compression case, the FEM simulation leads to poor results for coarse nodal discretization compared to MTLED using either CME or MMLS-SEBCIEM approximants (Figure \ref{fig:cube_compression}).

\subsection{Cylinder Under Locally Applied Indentation} \label{subsec:cylinder_indentation}

\begin{figure}[htbp]
    \centering
    \includegraphics[width=\textwidth]{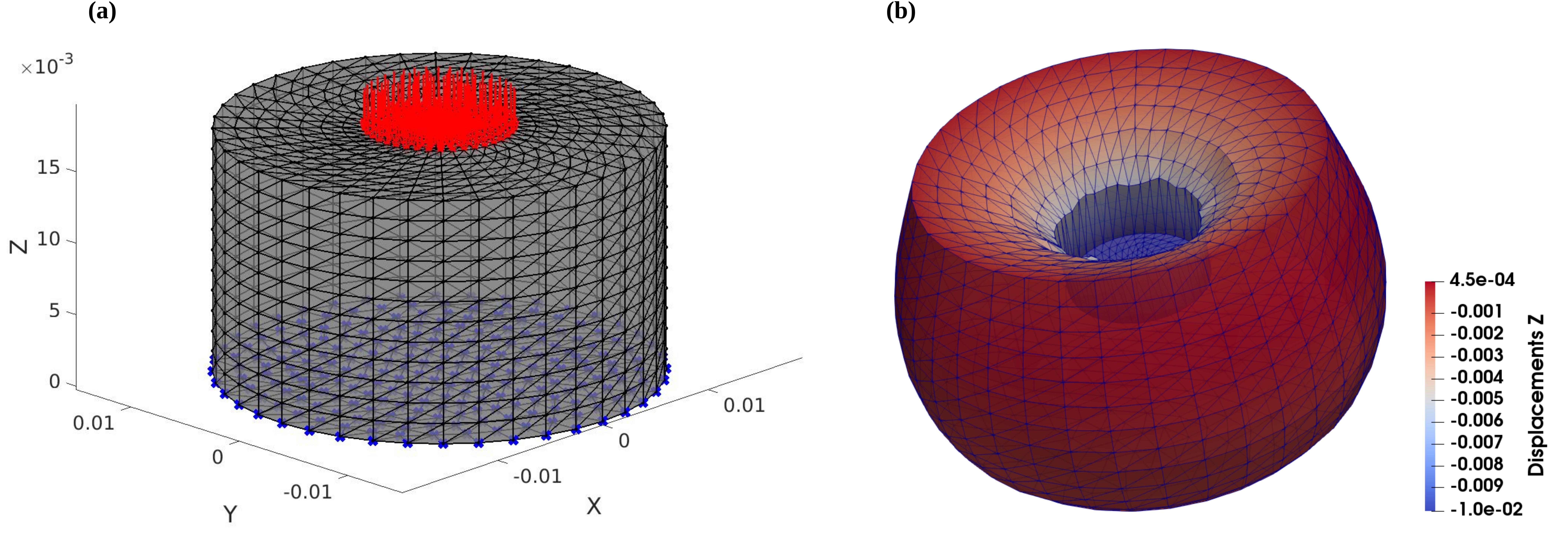}
    \caption{(a) Fixed (blue) and displaced (red) boundary nodes for $60\%$ locally applied indentation on a cylinder. (b) Steady state deformation of the cylinder under $60\%$ locally applied indentation. The colormap represents the resulting displacements at Z-axis direction for 177372 quadrature points.}
    \label{fig:cylinder_indent}
\end{figure}

We highlight the applicability of the proposed CME-MTLED method against extreme indentation and demonstrate its applicability beyond what is possible with FEM. In detail we consider indentation of a cylindrical domain with height $h=17 mm$ and diameter $d=30 mm$, modelled by the hyper-elastic neo-Hookean constitutive law as in \cite{Joldes2019}. The sub-region of the cylinder’s top surface (illustrated in Figure \ref{fig:cylinder_indent}) undergoes a vertical displacement of $10 mm$ ($60\%$ of the initial height) while the nodes belonging to the bottom surface are fixed to zero displacements at all directions. The assigned material properties are similar to the test cases in Subsections \ref{subsec:unconstrained_compression} and \ref{subsec:uniaxial_deformation} ($YM = 3000 Pa$; $v = 0.49$; $\rho=1000 kg/m^3$). The cylindrical domain consists of $8223$ nodes and $177372$ quadrature points (4-point quadrature rule). The CME approximants are constructed using the set-up of $N_R = 2$ and $s=2$; $m=3$; $a=2.0$; $k=2$. The deformation of the indented cylinder at steady state is given in Figure \ref{fig:cylinder_indent}. The described displacements at the bottom and top surfaces are exactly satisfied.

\begin{figure}[htbp]
    \centering
    \includegraphics[width=\textwidth]{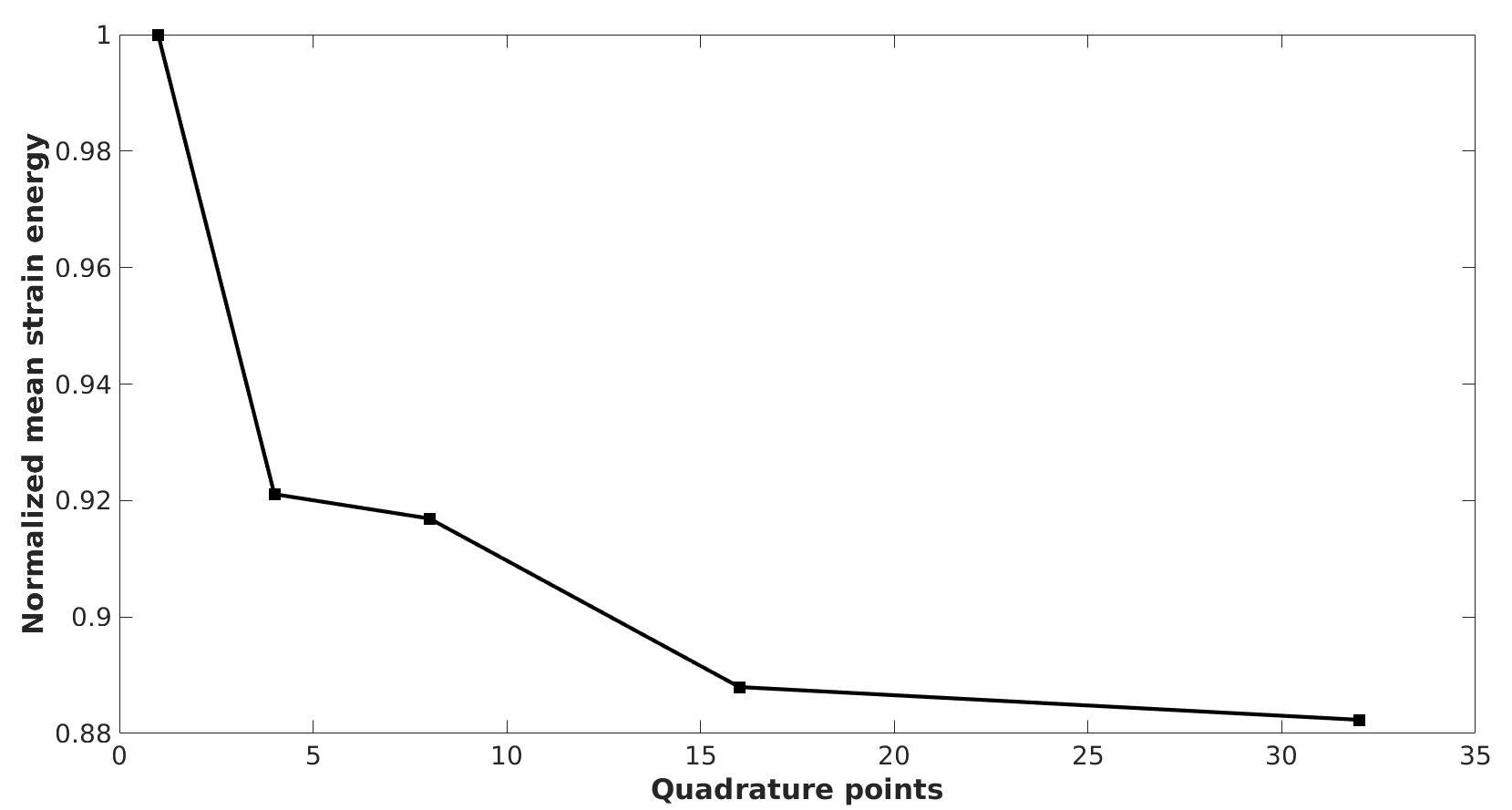}
    \caption{Mean values of strain energy for increasing quadrature of a locally applied indentation to a cylinder. The strain energy reduces monotonically for increasing quadrature points. Mean strain energy values are normalized to the maximum value (1 quadrature point) to aid the results’ interpretation.}
    \label{fig:cylinder_quad}
\end{figure}

The effect of the quadrature on the accuracy of the proposed method is evaluated. The mean strain energy density is measured for several simulations using 1, 4, 8, 16, and 32 quadrature points per integration cell. Quadrature with 8, 16, and 32 quadrature points is performed by using the 4-point quadrature rule after subdividing the cells of the original background mesh using 2, 4, and 8 subdivisions respectively. It is shown that for increasing the number of the quadrature points, the strain energy density is reduced monotonically (Figure \ref{fig:cylinder_quad}). The difference of the mean strain energy density for quadrature using 4 and 32 points per integration cell is found $3.9\%$. Therefore, performing integration with the 4-point quadrature rule appears to be a good compromise between accuracy and efficiency. Additional improvement of the accuracy with minimum computational overhead could be achieved by using the adaptive integration method described in \cite{Joldes2015b}.

\subsection{Cardiac Model Geometry Under Locally Applied Indentation} \label{subsec:cardiac_indentation}

\begin{figure}[htbp]
    \centering
    \includegraphics[width=\textwidth]{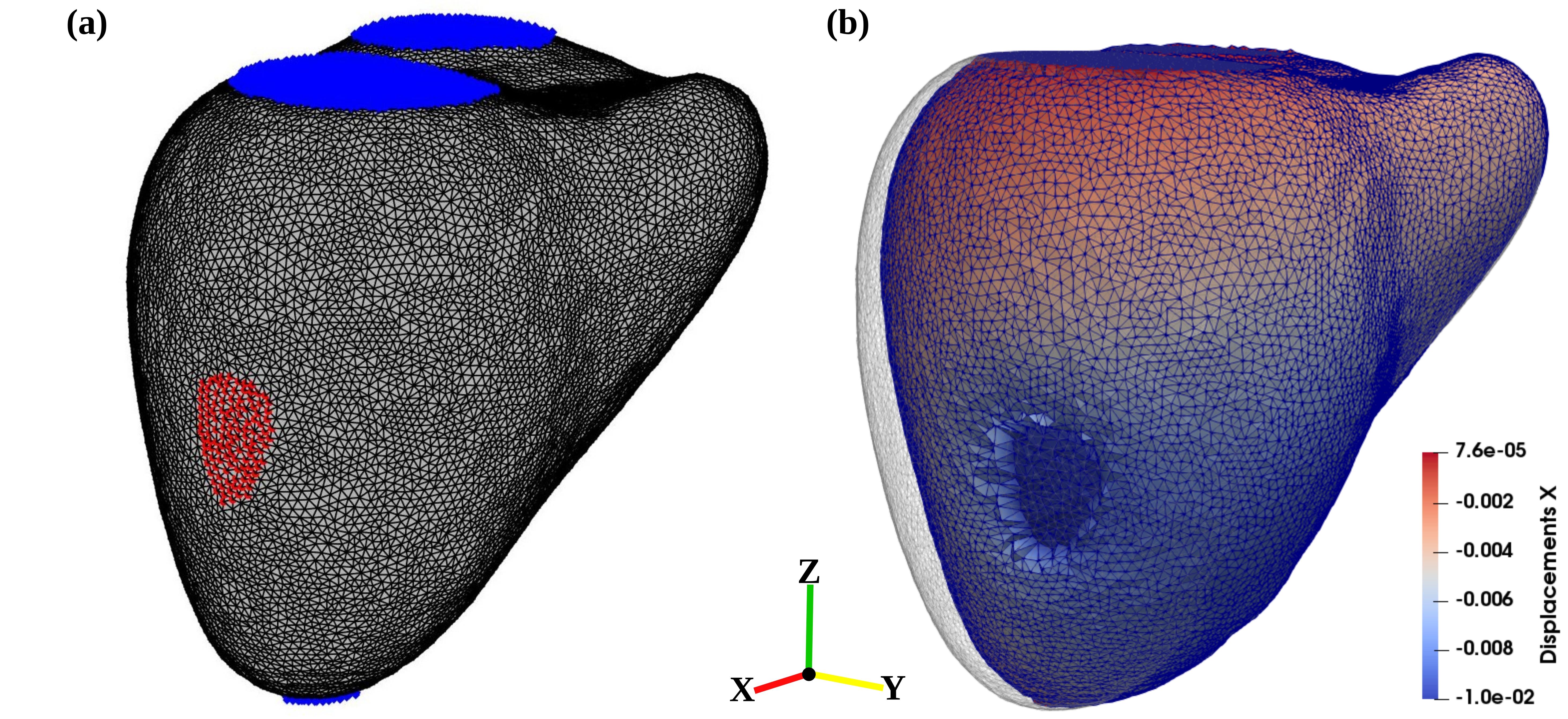}
    \caption{(a) Nodes (blue) at the top and bottom faces of the cardiac bi-ventricular model are fixed in all directions. The nodes on a small patch located at the right ventricle (red) are displaced by $u_x=0.01m$. (b) Steady state deformation during right ventricle indentation computed with the CME-MTLED using ($s = 2, a = 2.0$).}
    \label{fig:heart_indent}
\end{figure}

Changes in myocardium stiffness are related to different cardiac conditions or diseases (e.g., diastolic dysfunction \cite{Jackson2000}, tachycardia-induced cardiomyopathy \cite{Khasnis2005}). Evaluation of myocardial material properties under such circumstances is relevant for disease assessment and surgical planning. Material properties evaluation is based on the estimation of the force/displacement relation for a predefined displacement \cite{Shishido1998,Blair1996}. In this context, we simulate the indentation of a 3D cardiac bi-ventricular geometry. Our objective is to demonstrate the applicability of the CME-MTLED method on domains with arbitrary shape and its ability to satisfy the prescribed displacement conditions on such domains. We use a mesh composed of $38764$ nodes and $186561$ tetrahedral cells which was generated from MRI cardiac segmentation images using TetGen \cite{Si2015}. The segmentation images were made available from \cite{Bai2015}. Tissue indentation in the $x$-axis direction is simulated imposing a $u_x=0.01 m$ displacement on a patch region at the right cardiac ventricle while the bottom and top faces of the ventricles are constrained in all directions (Figure \ref{fig:heart_indent}). We use the neo-Hookean constitutive model with parameters $YM = 3000 Pa$; $v = 0.49$; $\rho=1000 kg/m^3$ to describe the mechanical response of the cardiac tissue. We assume that the cardiac bi-ventricular geometry is not beating. It should be noted that the selected essential boundary conditions and constitutive model are not based on available experimental data and do not represent neither in-vivo boundary conditions nor the realistic response of the heart. They are rather chosen aiming to demonstrate the capabilities of the CME-MTLED method. A realistic biomechanical simulation is out-of-scope of the current study.

The indentation simulation is performed using both the CME ($s=2$, $a=2.0$) and the MMLS-SEBCIEM. The characteristics $\Delta t_{crit}$, $N_{exe}$, $t_{exe}$ and the EBC values at steady state for the two approximation schemes are given in Table \ref{tab:heart_indent_results}. It is shown that using CME, the EBC can be imposed exactly without the need of any special treatment such as SEBCIEM or EBCIEM. In addition, solution stability can be achieved using longer explicit integration time steps leading to a $27\%$ reduction in the computational time.

\begin{table}[H]
    \caption{Simulation characteristics (critical explicit integration time step - $\Delta t_{crit}$, number of execution steps - $N_{exe}$, execution time - $t_{exe}$, Essential boundary conditions (EBC) values) for cardiac bi-ventricular model in indentation.}
    \centering
    \begin{threeparttable}
        \begin{tabular}{c c c c c c}
            \toprule
            Approx.  & $\Delta t_{crit}$ (ms) & $N_{exe}$ & $t_{exe}$ (h)  &  Displaced EBC\tnote{$\dagger$} (m) & Fixed EBC\tnote{$\dagger$} (m) \\
            \midrule
            CME &  0.636 & 36073 & 1.06 & -0.01 $\pm$ 1.0E-17 & 0.00 $\pm$ 0.00  \\
            MMLS-SEBCIEM & 0.455 & 49692 & 1.53 & -0.01 $\pm$ 3.3E-15 & 0.00 $\pm$ 1.7E-15 \\
            \bottomrule
        \end{tabular}
        \begin{tablenotes}
            \item [$\dagger$] Mean EBC values at steady state for fixed and displaced nodes with mean variation
        \end{tablenotes}
    \end{threeparttable}
    \label{tab:heart_indent_results}
\end{table}

\section{Concluding remarks} \label{sec:conclusions}

In the present study, we presented the extension of Cell-based Maximum Entropy (CME) approximants in the three-dimensional Euclidean space ($E^3$) using the R-functions technology based on the work of Millán et al \cite{Millan2015}, and Biswas and Shapiro \cite{Biswas2004}. The CME approximants were used specifically in the MTLED method, however they could be easily introduced in other solvers (explicit, implicit, semi-implicit). We evaluated the accuracy and efficiency of the CME in the Meshless Total Lagrangian Explicit Dynamics (MTLED) method through comparison with Modified Moving Least Squares (MMLS) and Finite Element Method (FEM). We presented two numerical examples to verify the capability of CME MTLED to generate accurate solutions to nonlinear equations of solid mechanics governing the behavior of soft, deformable tissues. We also verified the accuracy of the proposed scheme against extreme indentation, with the indentation depth reaching $60\%$ of the initial height of the sample. The evaluation was performed by analyzing the signed relative error in strain energy density for several combinations of parameters used in CME. Finally, we verified the accuracy of essential boundary conditions (EBC) imposition on domains of arbitrary shape in an indentation simulation of a cardiac bi-ventricular model derived from cardiac MRI segmentation images.

In all the numerical examples, the use of the CME approximants allowed longer explicit time integration step without compromising the numerical stability. A gain in computational time up to $27\%$ was achieved for CME, while demonstrating similar convergence to MMLS. The weak Kronecker-delta property of CME allowed to directly impose EBC, avoiding special EBC imposition treatments, which have been previously proposed to deal with the lack of the Kronecker delta property in MMLS. The smooth derivatives of CME led to derivatives matrix with smaller eigenvalues compared to MMLS and hence longer explicit time integration steps. CME was proven a valuable alternative to the group of MLS approximants in the context of MTLED and other similar EFG frameworks. A naive implementation of the CME approximants in $E^3$ is available to download at \texttt{\url{https://www.mountris.org/software/mlab/cme}}.

\section*{Acknowledgments}
This work was supported by the European Research Council through grant ERC-2014-StG 638284, by MINECO (Spain) through project DPI2016-75458-R and by European Social Fund (EU) and Arag\'on Government through BSICoS group (T39\_17R). Computations were performed by the ICTS NANBIOSIS (HPC Unit at University of Zaragoza). The authors wish to thank CONICET and grant PICTO-2016-0054 from UNCuyo-ANPCyT for funding the third author and the Department of Health, Western Australia, for its financial support to the fourth author through a Merit Award. The fifth and last authors acknowledge the funding from the Australian Government through the Australian Research Council ARC (Discovery Project Grants DP160100714, DP1092893, and DP120100402).

\bibliographystyle{unsrt}  
\bibliography{references}

\end{document}